\documentclass[12pt]{amsart}

\usepackage{latexsym}
\usepackage{psfrag}
\usepackage{amsmath}
\usepackage{amssymb}
\usepackage{epsfig}
\usepackage{amsfonts}

\newcommand{\R}{{\mathbb R}}

\newcommand{\kb}{\overline{k}}

\newcommand{\be}{\begin{equation}}
\newcommand{\ee}{\end{equation}}

\def\rr{\mathbb R}

\def\Si{\Sigma}
\def\la{\lambda}

\def\ssu{\subset}

\def\<{\langle}
\def\>{\rangle}

\def\0{{\mathbf 0}}

\def\.{\hskip.06cm}
\def\ts{\hskip.03cm}

\def\vol{{\text {\rm vol}}}

\def\pt{\partial}

\newtheorem{prop}{Proposition}[section]

\newtheorem{lemma}[prop]{Lemma}
\newtheorem{thm}[prop]{Theorem}

\newtheorem{remark}[prop]{Remark}

\newtheorem{claim}[prop]{Claim}

\newcommand{\II}{I\hspace{-0.1cm}I}

\newcommand{\kd}{\dot{k}}
\newcommand{\deltad}{\dot{\delta}}

\begin{document}

\title{Profiles of inflated surfaces}
\author{Igor Pak}
\address{School of Mathematics\\
University of Minnesota \\
Minneapolis, MN 55455\\
U.S.A.}
\email{pak@umn.edu}
\author{Jean-Marc Schlenker}
\address{Institut de Math\'ematiques, 
Universit\'e Toulouse III\\
31062 Toulouse cedex 9, France}
\email{schlenker@math.univ-toulouse.fr}
\date{July 28, 2009}

\begin{abstract}
We study the shape of inflated surfaces introduced
in~\cite{B1} and~\cite{P1}.  More precisely, we analyze
profiles of surfaces obtained by inflating a convex polyhedron,
or more generally an almost everywhere flat surface, with a
symmetry plane.  We show that such profiles are in a
one-parameter family of curves which we describe explicitly
as the solutions of a certain differential equation.
\end{abstract}

\maketitle

\section{Introduction} \label{s:intro}

The classical isoperimetric problem says that the sphere
is the unique surface of a given area which encloses the
maximal volume (see e.g.~\cite{BZ-book}).  However, when
the area is substituted with other functionals of the
intrinsic metric, the uniqueness disappears and the problem
becomes very difficult and sometimes intractable.

In this paper we study the geometry of \emph{inflated surfaces},
defined as the surfaces~$\Si$ of maximal volume among all
embedded surfaces homeomorphic to a sphere whose intrinsic
metric is a submetric to the intrinsic metric of~$\Si$.
In other words, we require that for every surface~$\Si'$
homeomorphic to~$\Si$ if all geodesic distances between
points in~$\Si'$ are less or equal to the geodesic distances
between the corresponding points in~$\Si$, then the volume
of~$\Si'$ is less or equal to that of~$\Si$.
This notion was originally introduced
by Bleecker~\cite{B1}, and further studied in~\cite{P1}.

Of course, every sphere is an inflated surface, since the
smaller geodesic distances imply smaller surface area.
However, not every convex surface is an inflated surface;
in fact, very few of them are.  For example, it was shown
in~\cite{B2} that not all ellipsoids with distinct axes
are inflated surfaces, and possibly none of them are.
To see an example of this phenomenon, consider a nearly
flat ellipsoid~$\Si$, with two large equal axes and one
very small one.  Think of the surface of~$\Si$ as of
non-stretchable closed balloon and blow some air into it.
Although the resulting surface is non-convex, it is still
isometric to~$\Si$ and encloses a larger volume, implying
that~$\Si$ is not an inflated surface.

The idea behind inflated surfaces is to start with a given surface~$S$
and blow air into it until no longer possible.  The resulting
surface~$\Si=\Si(S)$ will be the inflated surface, and will
typically be non-convex and not everywhere smooth even if~$S$
was convex and smooth.  Also, in the limit the geodesic distances
can (and typically do) become smaller in~$\Si$ than in~$S$,
so the inflated surface~$\Si$ (perhaps, counterintuitively) has
larger volume and smaller area than~$S$.
See sections~\ref{s:examples},~\ref{s:final}
and~\cite{P1} for many examples and conjectures on inflated surfaces.

Unfortunately, little is known about uniqueness and regularity of
inflated surfaces, even in the most simple cases.  We are motivated
by~\cite{Pa} (see also~\cite{MO}), where the shape of a
\emph{mylar balloon} was computed, defined as the inflated surface
of a doubly covered disc (glued along the boundary circle).  In this
case the symmetry can be utilized to obtain a complete profile by
solving a one-dimensional variational problem.  We consider a large
class of almost everywhere flat convex surfaces with a plane symmetry.
Examples include doubly covered regular polygons, cubes (or more
general brick surfaces), other Platonic and Archimedean solids, etc.
We then use tools from classical differential geometry and certain
heuristic assumptions to compute the \emph{profile} of the corresponding
inflated surfaces, defined as the shape of the curve in the symmetry plane.
We show that there is essentially one parametric family of such profiles,
which are all solutions of a special third order differential equation.
Given the variety of examples which include the mylar balloon, this
result may again seem counterintuitive.  This is the first result on
the shape of general inflated surfaces.

The rest of this paper is structured as follows.  In the next
section we give formal definitions and state the main theorem.
In Section~\ref{s:basic} we present a heuristic argument on
the metric of the inflated surfaces.\footnote{\, As we mention
later on, this argument can be made rigorous under certain regularity
assumptions.}  We use these results in Section~\ref{s:equations}
to write explicit equations on
the curvature of inflated surfaces obtained by an inflation of a
nearly flat surfaces.  In Section~\ref{s:symmetry} we use the
plane symmetry to prove the main result.  In the following
section (Section~\ref{s:examples}) we discuss a number of special cases
and make a special emphasis on the regularity assumptions.  We conclude
with final remarks in Section~\ref{s:final}.

\bigskip

\section{Definitions and main results} \label{s:main}

Let $S$ be a closed compact surface embedded in~$\rr^3$ and
homeomorphic to a sphere.  Throughout the paper we assume
our surfaces~$S$ are $C^5$ smooth everywhere except on a finite
union of curves in~$S$.  We further discuss the regularity
assumptions in Section~\ref{s:final}.

We say that the surface $S'\ssu \rr^3$ is \emph{submetric} to~$S$,
write $S' \preccurlyeq S$, if there exists a Lipschitz~1
homeomorphism \. $f: S \to S'$, i.e. such that the geodesic
distance satisfies \.
$|f(x)f(y)|_{S'} \le |xy|_{S}$ for all $x,y \in S$.
In particular, if~$S'$ is isometric to~$S$, then it is also
submetric.  Define
$$\upsilon(S) \. = \, \sup_{S' \preccurlyeq S} \. \vol(S'),
$$
where here and throughout the paper by $\vol(S)$ we denote
the volume enclosed by the surface~$S$.

It is natural to assume, and has been explicitly conjectured
in~\cite{P1}, that when~$S$ is convex there is a
unique (up to rigid motions) embedded surface~$\Si = \Si(S)$ which
attains the supremum: \ts $\vol(\Si) = \upsilon(S)$.
The surface $\Si$ is the \emph{inflated surface}, and we also
refer to~$\Si$ as the \emph{inflation of~$S$}.

From now on we consider only surfaces~$S$ which are convex and
almost everywhere flat.  The example include the surfaces $\pt P$ of
convex polytopes $P \ssu \rr^3$, doubly covered convex plane
regions.  Our goal is the description of~$\Si = \Si(S)$, which we
assume to be uniquely defined and satisfy the above regularity assumptions.

We are now ready to present the main result of this paper.
Suppose~$S$ is symmetric with respect to a plane~$H$, and let
$C=\Sigma\cap H$ be the \emph{profile} of~$\Sigma$.  We assume
that~$\Sigma$ is $C^5$ smooth in the neighborhood of~$C$.
(if~$C$ contains finitely many non-smooth points, consider
a portion of~$C$ between them).
Let $k=k(t)$ be the geodesic curvature of~$C$, considered
as a curve in~$H$, and parameterized by the length of~$C$.

\begin{thm} \label{t:main}
The curvature~$k(t)$ of the profile~$C$ satisfies
the following differential equation:
\begin{equation}
  \label{eq:diff}
 k(t) \ts k'''(t) \. - \. k'(t) \ts k''(t) \. + \.
 k^3(t) \ts k'(t) \, = \. 0\..
\end{equation}
\end{thm}

As a corollary, we conclude that non-constant solutions of
$k(t)$ for which $k(0)=0$ are given by the following integral formula:
$$
\int_0^{k(t)} \frac{du}{\sqrt{(\mu-\lambda^2)+4\lambda u^2 - u^4}}
\, = \, \pm \. \frac{t}{2}\,,
$$
where $\lambda,\mu\in \R$ are constants.

Note that the solutions are invariant under the
change $k\to \lambda \ts k$,  $t\to t/\lambda$.
This can be seen immediately from the invariance of
equation~(\ref{eq:diff}) under
the same change of variables.  Therefore, there is only a
one-parameter family of possible profiles, up to dilation.
The examples and other applications of the theorem will
be given in Section~\ref{s:examples}.

\bigskip

\section{Basic description}  \label{s:basic}

Let $S\ssu\rr^3$ be an almost everywhere flat convex surface as
in the previous section.  For simplicity, the reader can always
assume that $S$ is the surface of a convex polytope, even though
our results hold in greater generality.  As before, we denote
by $\Si$ the inflation of~$S$.  We call~$g_0$ and~$g$ the metrics
on $S$ and~$\Si$, respectively.  We also call $J_0$ and $J$ the
conformal structures of the metrics $g_0$ and $g$.  By definition,
$g$~is obtained by contracting $g_0$ in some directions
(at some points) and has maximal interior volume under this condition.
Call~$\Omega$ the subset of $\Sigma$ where some direction is
contracted.  So~$\Omega$ is an open subset of~$\Sigma$.

\begin{claim}
  \begin{enumerate}
  \item At each point $x\in \Omega$ one direction is not contracted. We
call $\xi$ a unit vector in the non-contracted direction, and $\sigma$
the contraction factor in the direction orthogonal to $\xi$.
\item The integral curves of $\xi$ are geodesics for $g$.
\item They are also geodesics for $g_0$.
\item At each point of $\Omega$, $\II(\xi,\xi)\geq 0$.
  \end{enumerate}
\end{claim}

\begin{proof}
For~$(1)$, note that if all directions are contracted it is
possible to deform a little the surface $\Sigma$ by making a little
``bump'', then $g$ remains smaller than $g_0$ while the interior
volume increases.

Part~$(2)$ follows from the fact that if integral curves of~$\xi$
are not geodesics in the neighborhood of some point $x_0$, then it is
possible to deform $\Sigma$ by adding a little ``bump'' while
``correcting'' the variation of the integral curves of $\xi$
by deforming them towards their concave side on $\Sigma$. All
those curves then keep the same length and~$g$ remains smaller
than $g_0$ while the volume increases.

For~$(3)$, let~$\gamma$ be a segment of integral curve of $\xi$,
so that by~$(2)$ it is geodesic for~$g$.  After replacing~$\gamma$
by a shorter segment if necessary, we can assume that~$\gamma$
is the unique minimizing segment between its endpoints.
Let~$\gamma_0$ be a minimizing segment for $g_0$ between the
endpoints of $\gamma$. Denote by $L_g$ the length
for~$g$, and by $L_{g_0}$ the length for $g_0$. Then
$L_g(\gamma)\leq L_g(\gamma_0)$ since $\gamma$ is minimizing
for~$g$, while $L_{g_0}(\gamma_0)\leq L_{g_0}(\gamma)$ since
$\gamma_0$ is minimizing for $g_0$.

On the other hand, from above, metric~$g$ is obtained
from~$g_0$ by contracting some directions,
so that $L_g(\gamma_0)=L_{g_0}(\gamma_0)$.  Putting the above
inequalities together, we obtain:
$$
L_g(\gamma)\leq L_g(\gamma_0) \leq L_{g_0}(\gamma_0) \leq
L_{g_0}(\gamma)\,.
$$
But we also know that $L_g(\gamma)=L_{g_0}(\gamma)$, precisely
because $\gamma$ is an integral curve of $\xi$. It follows
from here that the three inequalities above are actually
equalities, so that $\gamma_0=\gamma$.  This finishes the
proof of part~$(3)$.

For~$(4)$, if $\II(\xi,\xi)<0$, then it is possible
to add a small ``bump'' increasing the volume and shortening
all integral curves of $\xi$, a contradiction.
\end{proof}

\begin{remark} {\rm
Let us note that if $S$ and $S'$ are the doubly covered polygons
$Q$ and~$Q'$, respectively, and $Q\subset Q'$, then the surface
obtained by inflating~$S$ is not necessarily contained in the surface
obtained by inflating~$S'$.  For example, take~$Q$ to be a square
and~$Q'$ to be a very thin rectangle of length almost the diagonal
of~$Q$.  On the other hand, the volumes of the inflated
surfaces are monotonic in this case: $\upsilon(S) \le \upsilon(S')$
(cf.~\cite[$\S 5.1$]{P1}).
}
\end{remark}

\bigskip

\section{Equations}  \label{s:equations}

In this section we present the equations satisfied by the data
describing~$g$ and~$\Sigma$ as a surface in $\R^3$.  The
basic hypothesis here is that $\Sigma$ is the image of~$S$
by a contracting immersion
(singular at the cone points of~$\Sigma_0$ of course)
which is a critical point of the volume among
such contracting immersions.

Note first that the integral curves of~$\xi$ are not necessarily
parallel on~$S$; this apparently happens a lot in interesting situations
but for the moment we stick to a more general setting. Consider such
a line~$\Delta_0$, and another such line~$\Delta$ very close to it.
Let~$y$ be the function on $\Delta_0$ defined as the distance between
$\Delta_0$ and~$\Delta$ along the normal to~$\Delta_0$, and set
$$ \rho = y'/y\,. $$
Then $\rho$ makes sense as a limit $y\rightarrow 0$.

To simplify notations we use a prime for derivation along $\xi$, i.e.,
for any function $f$ on $\Sigma$, $f'=df(\xi)$. We use a dot for the
derivation along $J_0\xi$, i.e., $\dot{f} = df(J_0\xi)$.

\begin{lemma}
The curvature of~$g$ is given by the following equation:
$$ K \, = \, - \ts \sigma''/\sigma \. - \. 2\ts \rho \ts \sigma'/\sigma\,.
$$
\end{lemma}

\begin{proof}
The distance between $\Delta_0$ and $\Delta$ (for~$g$ now)
is $\sigma \ts y$.  Since both lines are geodesics, a basic
argument on Jacobi fields on Riemannian surfaces gives
that $K \ts = \ts -(y\ts \sigma)''/(y\ts \sigma)$, and
the result follows from \ts $y''=0$.
\end{proof}

\medskip

It is necessary to consider the second fundamental form of
$\Sigma$, we write its coefficients in the basis $(\xi,J\xi)$ as
$\binom{k_0 \, \delta}{\delta \ k_1}$.

\begin{lemma} \,
$K \ts = \ts k_0\ts k_1 \ts - \ts \delta^2$.
\end{lemma}

\begin{proof}
This is the \emph{Gauss formula}.
\end{proof}

\medskip

\begin{lemma} The curvature~$k_0$ is proportional to $y\ts \sigma$ along
integral curves of~$\xi$.  In other words,
$$
k_0'/k_0 \, = \, \sigma'/\sigma \. + \. y'/y\,.
$$
\end{lemma}

\begin{proof}
Consider two bumps along an integral curve of $\xi$, in
opposite directions so that the total length is not changed.
The variation of length is proportional to~$k_0$ and to the
normal displacement, while the variation of volume is proportional
to $y\ts \sigma$ times the normal displacement.
This implies the result.
\end{proof}

\medskip

\begin{prop} Let $\nabla$ be the Levi-Civita connection of $g$. Then
\smallskip
{
$\nabla_\xi\xi = 0$, while
$\nabla_{J\xi}\xi = (\sigma'/\sigma +\rho) J\xi$.
}
\end{prop}

\begin{proof}
The first point is a direct consequence of point (3) of Claim 3.1, since
$\xi$ is a unit vector field for $g$.

For the second point,
let $\nabla^0$ be the Levi-Civita connection of $g_0$. The definition of $\rho$
shows that $\nabla^0_{J_0\xi}\xi = \rho J_0\xi$, while Claim 3.1 shows that
$\nabla^0 _\xi (J_0\xi)=0$. Therefore
$$ [\xi,J_0\xi] = \nabla^0_\xi(J_0\xi)-\nabla^0_{J_0\xi}\xi = -\rho J_0\xi~. $$
By definition of $\sigma$, $g(J_0\xi,J_0\xi)=\sigma^2$, so that $J_0\xi=\sigma J\xi$,
so that
$$ [\xi, J\xi] = \left[\xi,\frac{1}{\sigma} J_0\xi\right] =
\frac{1}{\sigma}[\xi, J_0\xi] + \left(\xi.\frac{1}{\sigma}\right)J_0\xi =
- \frac{\rho}{\sigma} J_0\xi -\frac{\xi.\sigma}{\sigma^2} J_0\xi =
-\rho J\xi - \frac{\sigma'}{\sigma} J\xi~, $$
and the result follows.
\end{proof}

\medskip

\begin{lemma}
The \emph{Codazzi equation} can be written as:
$$
\delta' \. - \. \kd_0/\sigma \. + \. 2\ts \delta \ts (\sigma'/\sigma \ts + \ts\rho) \, = \, 0\,, \ \text{and}
$$
$$ k_1' \. - \. \deltad/\sigma \. + \. (k_1-k_0)(\sigma'/\sigma \ts + \ts \rho) = 0\,.
$$
\end{lemma}

\begin{proof}
Let $B$ be the shape operator of $\Sigma$, i.e., $\II(u,v)=I(Bu,v)$.
Then the Codazzi equation is
$$
(\nabla_uB)(v) -  (\nabla_vB)(u) \, =\, 0\,,
$$
where $\nabla$ is the Levi-Civita connection on~$\Sigma$, and
$u,v$ are any two vector fields.  Writing this for $\xi, J\xi$
we get:
$$ \nabla_\xi (BJ\xi)  -  B(\nabla_\xi (J\xi)) - \nabla_{J\xi}(B\xi)
+ B(\nabla_{J\xi}\xi) \, = \, 0\,. $$
Now expressing this in terms of the coefficients of $\II$, we obtain:
$$ \nabla_\xi (\delta \xi + k_1 J\xi) \ts -
\nabla_{J\xi}(k_0\xi + \delta J\xi)
\ts + \. B((\sigma'/\sigma +\rho)J\xi) \, = \, 0\,.
$$
Therefore,
$$
\delta' \xi + k'_1 J\xi - dk_0(J\xi) \xi - k_0 (\sigma'/\sigma+\rho) J\xi  -
d\delta(J\xi) J\xi + \delta (\sigma'/\sigma +\rho) \xi
+ (\sigma'/\sigma +\rho)(\delta \xi + k_1J\xi) = 0\..
$$
On the other hand, $\sigma J\xi=J_0\xi$, so $d\delta(J\xi) = \deltad/\sigma$,
and similarly $dk_0(J\xi) = \kd_0/\sigma$. Separating the terms in~$\xi$
and in~$J\xi$ gives
$$
\delta' \. - \. \kd_0/\sigma \. + \. \delta \ts (\sigma'/\sigma \ts +\ts \rho)
\. + \. (\sigma'/\sigma \ts +\ts \rho)\ts \delta \. = \. 0\..
$$
We conclude:
$$k'_1 \. - \. k_0 \ts (\sigma'/\sigma\ts +\ts \rho) \. - \.
\deltad/\sigma \. + \.(\sigma'/\sigma\ts +\ts\rho)\ts k_1 \. = \. 0\.,
$$
which implies the result.
\end{proof}

\bigskip

\section{Lines of symmetry}  \label{s:symmetry}

In this section we consider a special case as in the main theorem
(Theorem~\ref{t:main}), when~$S$ has a symmetry plane~$H$.
Consider an integral curve of $\xi \ssu \Si \cap H$, which
means that it is a line of symmetry of the inflated surface.
For example, in a mylar balloon, all segments going through the
center are such lines of symmetry.  Other examples include the
case when~$S$ is the doubly covered regular polygon or a rectangle.
In each case plane~$H$ is orthogonal to the polygon and is a
symmetry plane of both~$S$ and~$\Si$.

We consider the case of a regular polygon for simplicity; as the
reader will see the general case follows by the same argument.
Let~$Q$ be a regular $n$-gon, made by gluing~$n$ copies of a
triangle~$T = (OAB)$. The copies of~$T$ are all glued so that
their vertices~$O$ are glued together. Moreover the triangle
is symmetric with respect to the line orthogonal to~$AB$
going through~$O$.  Denote by~$E$ be the midpoint of~$AB$.

As before, we assume that there is a unique surface~$\Sigma$
which has maximal volume and is submetric to~$S$.  It follows
from uniqueness that this surface has all the symmetries
of $\Sigma_0$. So it sufficient to study the quantities describing
the situation on~$T$.

The equations above simplify somewhat when considered on an axis
of symmetry, for instance on the segment $OE$ of the triangle
considered above. Then $\delta=0$, and $\kd_0=\kd_1=0$.
Such lines are lines of curvature (integral lines of the curvature directions),
we suppose that the corresponding principal curvature is~$k_0$.
In that case the three basic equations reduce to the following:
\begin{itemize}
\item \ for the Gauss equation,
\begin{equation} \label{eq:gauss}
k_0 \ts k_1 \. = \. -\ts \frac{(\sigma \ts y)''}{\sigma \ts y}\,,
\end{equation}
\item \ for the Codazzi equation,
\begin{equation} \label{eq:codazzi}
k'_1 \. = \. \left(\frac{\sigma'}{\sigma} \ts +
\ts \frac{y'}y\right) (k_0-k_1)\,,
\end{equation}
\item \ for the ``conservation of curvature'',
\begin{equation} \label{eq:conservation}
\frac{k'_0}{k_0} \. = \. \frac{\sigma'}\sigma \ts + \ts \frac{y'}y\,.
\end{equation}
\end{itemize}
One can use this last expression in the previous two to get
$$
k_0 \ts k_1 \. = \. - \ts \frac{(\sigma \ts y)''}{\sigma \ts y}
\. = \. - \left(\frac{(\sigma \ts y)'}
{\sigma\ts y}\right)' \ts -  \left(\frac{(\sigma \ts y)'}{\sigma\ts y}\right)^2
\. = \. -\left(\frac{k'_0}{k_0}\right)' \ts-\left(\frac{k'_0}{k_0}\right)^2
\. = \. -\ts\frac{k''_0}{k_0}\.,
$$
$$ k_0 \ts k'_1 \. = \. (k_0-k_1)\ts k'_0\..
$$
Replacing~$k_1$ in those
equations, we find after simple computations that~$k_0$ satisfies
equation~(\ref{eq:diff}), which we recall for convenience:
\begin{equation}\label{eq:z}
k \ts k''' \ts - k'\ts k''\ts +\ts k^3\ts k'\ts =\ts 0\..
\end{equation}

As we mentioned in Section~\ref{s:basic},
this equation is invariant under the transformation
$k \mapsto \lambda \ts k$, $t\mapsto t/\lambda$, which makes sense
since this homogeneity condition on~$k_2$ corresponds to the
invariance of the class of inflated surfaces under scaling.

The other key quantities describing the surface at the symmetry
line can then be recovered from~$k$. Setting $u=\sigma\ts y$, we get:
$$
\frac{u'}{u} \,= \,\frac{k'}k\,.
$$
From here we see that~$u$ is proportional to~$k$, while
$$ k_1 \, = \, -\ts\frac{k''}{k^2}\,.
$$

Finally, let us mention that equation~\ref{eq:z} can be solved
implicitly in the following cases.  The proof is straightforward.

\begin{prop}\label{p:eq}
The solutions of (\ref{eq:z}) vanishing at $t=0$ are $k=0$, and the
functions defined implicitly by
\begin{equation} \label{sol:0}
\int_0^{k(t)} \. \frac{ds}{\sqrt{(\mu-\lambda^2)\ts +\ts 4\ts\lambda\ts s^2 \ts -\ts s^4}} \, =
\, \pm \ts \frac{t}{2}\,,
\end{equation}
with constants $\lambda, \mu\in \R$. The solutions of~(\ref{eq:z}) with
$k(0)=k'(0)=0$ are $z=0$, and the functions defined by
\begin{equation} \label{sol:00}
\int_0^{k(t)} \frac{ds}{\sqrt{\lambda s^2 - s^4}} =
\pm \frac{t}{2}\,.
\end{equation}
\end{prop}

\begin{proof}
Let $k$ be a solution of~(\ref{eq:z}), defined on an interval $I\subset \R$.
Then
$$ \left(k\ts k'' -(k')^2 +\frac{k^2}4\right)'=0~, $$
so that
$$ k\ts k'' -(k')^2 +\frac{k^2}4=a $$
for some constant $a\in \R$.
Let $y$ be the inverse function of $k$ (on a subinterval $J$ of
$I$ where $k$ is monotonous), then $k'\circ y=1/y'$ and
$y'^2 k''\circ y + y'' k'\circ y=0$, so that
$$ -\frac{ty''(t)}{y'(t)^3} - \frac{1}{y'(t)^2} + \frac{t^4}{4} = a~. $$
Set $u(t)=ty'(t)$, then, on the sub-interval where $u(t)\neq 0$,
$$ \frac{u'(t)}{u(t)^3} = \left( \frac{t}{4}-\frac{a}{t^3} \right)$$
so that there exists $b\in \R$ such that
$$ \frac{1}{u(t)^2} = \left( \frac{a}{t^2} - \frac{t^2}{4}\right)+b~, $$
for some $b\in \R$, and
$$ u(t) = \pm \frac{2|t|}{\sqrt{4a-t^4+4bt^2}}~. $$
Therefore,
$$ y(t) = \pm \int^t \frac{2 ds}{\sqrt{4a-s^4+4bs^2}}~. $$
Since $y(k(r))=r$ by definition of $y$, $k$ satisfies the equation
$$ \int^{k(r)} \frac{2 ds}{\sqrt{4a-s^4+4bs^2}} = \pm\frac{r}{2}~, $$
which is equivalent to (\ref{sol:0}) by a simple change of the constants.

The solutions such that $k(0)=0$ correspond to taking the primitive
equal to $0$ at $0$, and those solutions have vanishing derivative at $0$
if and only if the integrand goes to $\infty$ at $s=0$, that is,
if and only if $a=0$.
\end{proof}

\medskip

\begin{remark} \label{r:curv} {\rm
Taking into account the invariance under homotheties,
the proposition implies that there is a one-parameter
family of possible ``profiles'', which we call~$C_\nu$.
For all the cases with planar symmetry obtained from a doubly covered
convex figure, we have $k'=0$ at the ``equator''.  In addition to this,
for the mylar balloon we have $k=0$ at the ``pole'', which determines
the curvature~$k$ uniquely.  For the other examples, however, there is
no reason to believe that $k=0$ at the ``pole'' so $k$ is completely really
determined, we still need one more boundary condition.
}
\end{remark}

\bigskip

\section{Examples}  \label{s:examples}

\subsection{The mylar balloon}

Consider a \emph{mylar balloon}, defined as above by gluing two
copies of a disk.  Both Theorem~\ref{t:main} and Proposition~\ref{p:eq}
can be applied in this case (see also Remark~\ref{r:curv}), and
this determines the balloon profile in terms of curvature (see below).
In fact, this is the only case when the profile was already
computed~\cite{MO,Pa}, but described in a different manner.

\begin{figure}[hbt]
\begin{center}
\epsfig{file=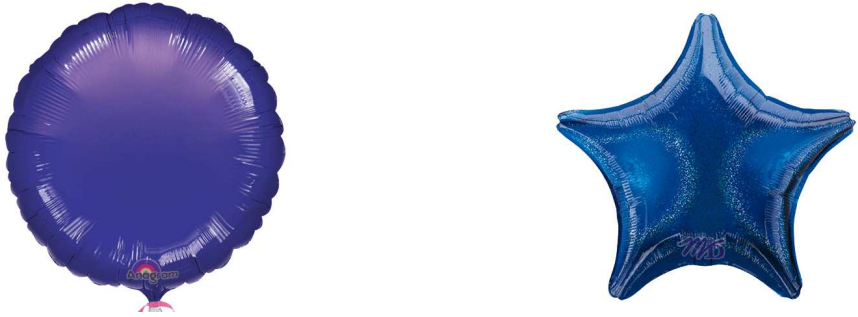,width=11.6cm}
\end{center}
\caption{Pictures of two party balloons.}
\label{f:balloons}
\end{figure}

Paulsen \cite{Pa} showed that the intersection with the upper right
quadrant of the profile of the mylar balloon
is the graph of the function $f:[0,a]\rightarrow \R$ given by
\begin{equation}
  \label{eq:paulsen}
 f(x) = \int_x^a \frac{u^2 du}{\sqrt{a^4-u^4}}~,
\end{equation}
where $a$ is the radius of the inflated balloon. This profile can be 
characterized by the following simple geometric property, which can be found 
(implicitly at least) in~\cite{MO}. We provide a direct proof here for the 
reader's convenience.

\begin{prop} \label{p:linear}
The profile of the mylar balloon is characterized by the fact that its
curvature is a linear function of~$x$: \ts $k(x)=-2x/a^2$, where~$a$ is the
radius of the balloon.
\end{prop}

\begin{proof}
Recall that the curvature of the graph of a function $f:\R\rightarrow \R$
is given by
$$ k(x) = \frac{f''(x)}{\bigl(1+f'(x)^2\bigr)^{3/2}}~. $$
Here
$$ f'(x) = -\frac{x^2}{\sqrt{a^4-x^4}}~, $$
so that
$$ \sqrt{1+f'(x)^2} = \frac{a^2}{\sqrt{a^4-x^4}}\,. $$
Similarly,
$$ f''(x) = -\frac{2x(a^4-x^4)+2x^3 \cdot x^2}{\bigl(a^4-x^4\bigr)^{3/2}} =
\frac{-2xa^4}{\bigl(a^4-x^4\bigr)^{3/2}}~, $$
and the result follows.
\end{proof}

With this proposition in mind it is easy to check that the curvature of the
profile of the mylar balloon is precisely the solution of (\ref{eq:z})
which vanishes at $t=0$ -- corresponding to the north pole of the balloon --
and with vanishing derivative at the equator, corresponding to $x=a$. Note
that we switch between a parameterization by the $x$ coordinate, in
(\ref{eq:paulsen}) and in Proposition~\ref{p:linear}, to a parameterization
by arclength, which is used in the description of the curvature as
solution of (\ref{eq:z}). For this reason we will denote by $\kb$ the
curvature of the profile as a function of $x$, keeping the notation
$k$ for the curvature as a function of the arclength parameter $t$.

We will now check that the profile in (\ref{eq:paulsen})
corresponds to the solution of (\ref{eq:z}) obtained by choosing
$\lambda=0,\mu=16/a^4$ in (\ref{sol:0}). For this we use
Proposition~\ref{p:linear}, which shows that
$$ \int_0^{\kb(x)} \frac{1}{\sqrt{16/a^4-s^4}} \. ds = \int_0^{-2x/a^2}
\frac{1}{\sqrt{16/a^4-s^4}} \. ds=
-\frac{1}{2} \int_0^x \frac{a^2}{\sqrt{a^4-u^4}} \. du~, $$
where the last equality uses the change of variables $s=-2u/a^2$.
Consider the arclength parameter $t$ on the profile, as a function
of $x$. Then
$$ t(x) = \int_0^x \sqrt{1+f'(x)^2} \. dx = \int_0^x \frac{a^2}{\sqrt{a^4-u^4}} \. du~.$$
As a consequence we obtain that
$$ \int_0^{\kb(x)} \frac{s^2}{\sqrt{16/a^4-s^4}} \. ds = -\frac{t(x)}{2}~, $$
so that the profile in (\ref{eq:paulsen}) corresponds to (\ref{sol:0})
for $\lambda=0,\mu=16/a^4$.

This profile has two prominent features.
\begin{itemize}
\item $k(0)=0$, this clearly follows from (\ref{eq:conservation}), since
$y'/y\rightarrow \infty$ at the north pole.
\item $k'(t)=0$ for $t$ corresponding to $x=a$, that is, to the equator of
the mylar balloon. This is clear since $x'(t)\rightarrow \infty$ as $x\rightarrow a$,
while $\kb'$ is bounded at that value of $x$ because $\kb$ is a linear function of $x$.
This condition should apply to the profile of the mylar balloon for a clear
symmetry reason.
\end{itemize}
Those two conditions characterize, up to dilation, the profile of the mylar balloon
among solutions of (\ref{eq:z}).

\subsection{The square and rectangular pillow}

A rectangular pillow is the surface obtained by inflating the doubly
covered rectangle. We call $v_\pm$ the centers of the two copies of the
rectangle, and $v_i, 1\leq i\leq 4$ its vertices. The rectangular pillow
is thus the union of eight triangles, of vertices $v_\pm, v_i, v_{i+1}$.
The square pillow is the special case of rectangular pillow for which
the doubly covered rectangle which is inflated is a square.
We make the following natural regularity hypothesis:  \emph{The
square pillow is $C^5$ on the interior of each of
the eight triangles.}

Clearly each of the triangles in the square pillow are congruent,
by the general hypothesis at the beginning of the section. We
consider one of them, say $(v_+,v_1,v_2)$. Let~$w$ be the
midpoint of $(v_1,v_2)$. $(v_+,w)$ is a line of symmetry of
the triangle, therefore a line of curvature, and the
corresponding principal curvature satisfies~(\ref{eq:z}).

To have a better understanding of this profile of the square pillow,
we can make a heuristically attractive hypothesis, which should
be compared with experimental data: we suppose that the surface is
$C^3$ at $w$. It then follows from symmetry that the derivative of $k$
at $w$ along $(v_+,w)$ is zero.

Under this hypothesis, the principal curvature $k$ along $(v_+,w)$
is a solution of (\ref{eq:z}) with vanishing derivative at $w$.
We conclude that the intersection
of the rectangular pillow with a plane of symmetry
containing $v_-,v_+$ but none of the $v_i$ is the union of two
copies of the profiles~$C_\nu$ for some~$\nu$, obtained from
equation~(\ref{sol:0}) by a simple integration once only one
additional parameter is known -- for instance the curvature of
the profile at $v_+$. If this curvature were zero at $v_+$, it
would imply that this curve is the same (up to dilation) as the
profile of the mylar balloon.

\begin{figure}[hbt]
\begin{center}
\epsfig{file=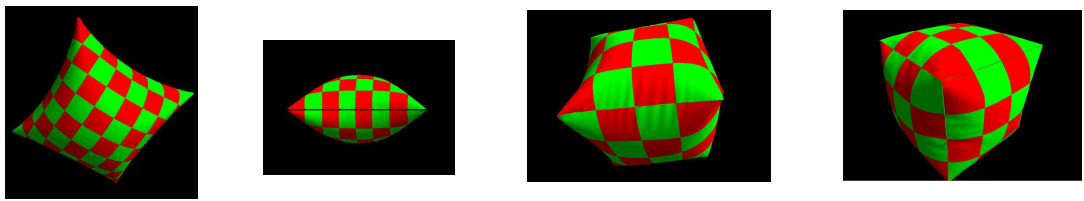,width=14.3cm}
\end{center}
\caption{Computer simulation of the square pillow and the inflated cube.}
\label{f:square-cube}
\end{figure}

Let us note that, to understand the more general case of rectangular
pillows, it would be useful to determine the dependence of
parameters $\la, \mu$ on the rectangle aspect ratio.  The case of
a square pillow is particularly attractive, and known as the
\emph{teabag problem} in recreational literature~\cite{K}.
Let us also mention the simulations by Gammel~\cite{G}
(see Figure~\ref{f:square-cube}),
and physical experiments by Robin for the conjectured formula for
the volume~\cite{R}.

\subsection{Doubles of polygons}

It is also possible to consider doubly covered regular $n$-gons
(cf.~Figure~\ref{f:balloons}).
We still call $v_\pm$ the centers of the two copies of the regular
polygon, and $v_i, 1\leq i\leq n$ their vertices (pairwise identified).
The inflated double $n$-polygon is the copy of $2n$ triangles, each
of which is $(v_\pm,v_i,v_{i+1})$. We make the same hypothesis as for
the square pillow, and obtain a profile shape which again only
depends on one parameter, for instance the curvature of the profile
at $v_\pm$, as for the square pillow (which is obviously a special case).

\subsection{The inflated cube}

Let us start now with the surface $S = \pt P$ or a unit cube~$P$.
Using the same analysis, each face of the inflated cube is cut into
four triangles and we can assume that each of them is $C^5$~smooth.
Under this hypothesis we find that the intersection of the inflated
cube with a plane of symmetry containing no vertex is the union
of four arcs of~$C_\nu$, for some~$\nu$.

\bigskip

\section{Final remarks}  \label{s:final}

\subsection{Why is this a proper mathematical model? }

There is a rather subtle point that needs to be made in favor
of the relevance of the model studied here for the physical
problem of understanding inflated surfaces.

It might appear at first sight that the model considered here,
based on contracting embeddings of a surface with maximal
volume, is quite different from what happens for true
inflated surfaces. Indeed for those surfaces there is no
contraction of the metric on the surface, but rather some
``plaids'' appear (see the wrinkles in Figure~\ref{f:balloons}).
One feature of those plaids is that they are necessarily
along geodesics on the surface, and this seems to impose
a constraint not present in the mathematical model.
However, there is a very good fit between the mathematical
models and the inflated surfaces as they are observed.

We believe that the resolution of this apparent paradox
lies in parts $(2)$ and~$(3)$ of Claim~3.1.  These
parts state precisely that the integral curves of the
non-contracted directions are geodesics on the surface.
In other terms, there is a mathematical constraint,
coming from the maximality of the surface, which happens
to coincide precisely with the physical constraints that
the plaids have to be along lines. For this reason the
``mathematical'' inflated surfaces are very close to the
observed ones, in spite of apparently different constraints
on their geometry.

\subsection{Future directions }

Perhaps the main open problem is to show that the inflated
surfaces are well defined and uniquely determined (see~\cite{P1}
for a complete statement).  Presumably,
this would imply the symmetry assumptions we made throughout
the paper.  Unfortunately, even in the case of the Mylar balloon
or rectangular pillow this is completely open.  It would also be
important to prove the regularity conditions, in particular
a formal proof of Claim~3.1.

Even under the uniqueness and regularity conditions, this paper
goes only so far towards understanding of the inflated pillow shapes.
Although we do not wish to suggest that in the case of rectangular
pillows the shape of the surface can be expressed by means of
classical functions, as in the case of the Mylar balloon~\cite{MO},
it is perhaps possible that it is a solution of an elegant problem
which completely defines it.  For example, the linearity of the
curvature as in Proposition~\ref{p:linear} completely characterizes
the Mylar balloon.  It would be interesting if the shape of rectangular
pillows has a similar characterization.

Finally, the \emph{crimping density} and \emph{crimping ratio}
(the ratio of the
area of inflated over non-inflated surface defined in~\cite{P1})
are interesting notions with potential applications to  material
science.  Exploring them in the case of rectangular pillows would
be of great interest.


\vskip.5cm

\noindent
{\bf Acknowledgements.} \. We would like to Andreas Gammel who
kindly allowed us to use his computer simulation figures, and to
Shelly Jenson from \ts {\tt BalloonManiacs.com} \ts for the permission
to use their high quality pictures of mylar balloons.
The first author was partially supported by the NSF and the NSA.
Both authors were supported by the MIT--France Seed Fund which
allowed this collaboration to begin.

\vskip.9cm




\end{document}